\documentclass[11pt]{amsart}
\usepackage{amssymb}
\usepackage{amsmath}
\usepackage[hidelinks]{hyperref}
\usepackage{etoolbox}
\usepackage{enumitem}
\usepackage{xcolor}
\usepackage{tikz-cd}
\usepackage[normalem]{ulem}
\usepackage{fancyvrb}

\newtheorem{theorem}{Theorem}[section]
\newtheorem{lemma}[theorem]{Lemma}

\newtheorem{proposition}[theorem]{Proposition}

\theoremstyle{remark}
\newtheorem{remark}[theorem]{Remark}

\newcommand{\qp}{\vskip .12cm}
\newcommand{\hp}{\vskip .2cm}
\newcommand{\p}{\vskip .4cm}

\newcommand{\F}{\mathbb{F}}
\newcommand{\Z}{\mathbb{Z}}
\newcommand{\Q}{\mathbb{Q}}
\newcommand{\R}{\mathbb{R}}
\newcommand{\C}{\mathbb{C}}
\newcommand{\Cc}{\mathbb{C}^{\times}}

\newcommand{\CO}{\mathcal{O}}

\newcommand{\til}{\tilde}
\newcommand{\bsl}{\backslash}
\newcommand{\ra}{\rightarrow}
\newcommand{\ira}{\hookrightarrow}
\newcommand{\sra}{\twoheadrightarrow}
\newcommand{\xra}{\xrightarrow}

\newcommand{\Lie}{\operatorname{Lie}}
\newcommand{\rank}{\operatorname{rank}}

\newcommand{\Ad}{\operatorname{Ad}}

\newcommand{\Spec}{\operatorname{Spec}}
\newcommand{\Frob}{\mathrm{Frob}}
\newcommand{\Irr}{\operatorname{Irr}}
\newcommand{\Tr}{\operatorname{Tr}}

\newcommand{\Id}{\operatorname{Id}}
\newcommand{\Lg}{\mathfrak{g}}
\newcommand{\Lh}{\mathfrak{h}}

\newcommand{\se}{\mathsf{e}}

\newcommand{\matr}[1]{\left[\begin{matrix}#1\end{matrix}\right]}

\newcommand{\Sym}{\operatorname{Sym}}
\newcommand{\supp}{\operatorname{supp}}
\newcommand{\Stab}{\operatorname{Stab}}

\begin{document}

\title{Geometric wave-front set may not be a singleton}
\author{Cheng-Chiang Tsai}
\date{Dedicated to Benedict H. Gross}
\thanks{The author is supported by MOST grant 110-2115-M-001-002-MY3.}
\email{chchtsai@gate.sinica.edu.tw}
\address{Academia Sinica, Institute of Mathematics, 6F, Astronomy-Mathematics Building, No. 1, Sec. 4, Roosevelt Road, Taipei 106319, Taiwan, and\vskip.2cm
National Sun Yat-Sen University, Department of Applied Mathematics, No. 70, Lienhai Rd., Kaohsiung 80424, Taiwan}
\begin{abstract}
	We show that the geometric wave-front set of specific half-integral-depth supercuspidal representations of ramified $p$-adic unitary groups is not a singleton.
\end{abstract}
\makeatletter
\patchcmd{\@maketitle}
{\ifx\@empty\@dedicatory}
{\ifx\@empty\@date \else {\vskip3ex \centering\footnotesize\@date\par\vskip1ex}\fi
	\ifx\@empty\@dedicatory}
{}{}
\patchcmd{\@adminfootnotes}
{\ifx\@empty\@date\else \@footnotetext{\@setdate}\fi}
{}{}{}
\makeatother

\maketitle 

\vskip-1cm $\;$

\tableofcontents

\vskip-1cm $\;$

\section{Introduction}

Let $F$ be a finite extension of $\Q_p$ and $G$ be a connected reductive group over $F$. For an irreducible smooth $\C$-representation $\pi$ of $G(F)$, the local character expansion of Howe and Harish-Chandra \cite[Thm. 16.2]{HC} asserts that the character $\Theta_{\pi}$ enjoys an asymptotic expansion on some neighborhood $U$ of the identity. To be precise, there exist constants $c_{\CO}(\pi)\in\C$ indexed by nilpotent $\Ad(G(F))$-orbits $\CO\subset\Lie G(F)$ such that
\begin{equation}\label{LCE}
\Theta_{\pi}|_U=\sum_{\CO}c_{\CO}(\pi)\cdot\left(\hat{I}_{\CO}\circ\log|_U\right)
\end{equation}
where $I_{\CO}$ is the orbital integral on $\CO$ and $\hat{I}_{\CO}$ its Fourier transform. Here we fix an $\Ad(G)$-equivariant isomorphism between $\Lg:=\Lie G$ and its dual. 

In \cite{MW87}, M\oe glin and Waldspurger generalized a result of Rodier \cite{Rod75} and showed that if $\CO$ is maximal among those with $c_{\CO}(\pi)\not=0$, then $c_{\CO}(\pi)$ is the dimension of the degenerate Whittaker model for $\pi$, in particular a positive integer. The set of those $\CO$ with $c_{\CO}(\pi)\not=0$ and maximal among such is commonly called the wave-front set, which has been the subject of many studies and results. For instance, \cite{Moe96} showed that for $p$-adic classical groups, any member of a wave-front set is {\it special} in the sense of Lusztig \cite{Lus79}. See also \cite{BM97}, \cite{JLS16}, \cite{Wal18}, \cite{Wal20}, \cite{GGS21}, \cite{Oka21a}, \cite{CMBO21}, \cite{AGS22a}, \cite{CMBO22a}, \cite{JLZ22a} and many more, particularly \cite{GGS21} for global discussions. Among the thread, there has been the conjecture that the wave-front set is contained in a single $\Ad(G(F^{sep}))$-orbit, or that the ``geometric wave-front set'' is a singleton. We give a counterexample to this conjecture.

Let $p=3$ and $F=\Q_3$ be our $p$-adic field, $E/F$ any ramified quadratic extension, and $k:=\F_3$ the common residue field of $E$ and $F$. Let $G=U_7(E/F)$ be the ramified unitary group over $F$ that splits over $E$. Explicitly, we identify $G(F)$ as the group of unitary operators on $E^7$ using the hermitian inner product
\begin{equation}\label{form1}
\langle x,y\rangle=x_1\bar{y}_7+...+x_7\bar{y}_1.
\end{equation}
Here, for $y\in E$ we denote by $\bar{y}$ its conjugate over $F$. Denote by $\CO_E\subset E$ the ring of integers in $E$ and $\mathfrak{m}_E$ the maximal ideal. Consider the filtration $(G(F)_r)_{r\in\frac{1}{2}\Z_{\ge 0}}$ given by
\[
G(F)_r:=\{g\in G(F)\subset M_{7\times 7}(E)\;|\;g-\Id_7\text{ has entries in }\mathfrak{m}_E^{2r}\}.
\]
This is the Moy-Prasad filtration at a specific special vertex (apart from a difference of index $2$ for $G(F)_0$, which will not concern us). Every subgroup within the filtration is normal in any previous group. Let us fix $\varpi\in E$ a uniformizer with $\bar{\varpi}=-\varpi$ once and for all. For any $n\in\Z_{\ge 0}$, the map
\begin{equation}\label{map1}
G(F)_{n+\frac{1}{2}}/G(F)_{n+1}\ira (\Id_7+\mathfrak{m}_E^{2n+1}M_{7\times 7}(E))/(\Id_7+\mathfrak{m}_E^{2n+2}M_{7\times 7}(E))\cong M_{7\times 7}(k)
\end{equation}
given by dividing the entries by $\varpi^{2n+1}$ is an injective map that identifies $G(F)_{n+\frac{1}{2}}/G(F)_{n+1}$ with the abelian group of $7\times 7$ self-adjoint matrices over $k=\F_3$. Here the adjoint is defined with respect to the form on $k^7$ induced by (\ref{form1}), namely
\begin{equation}\label{form2}
	\langle x,y\rangle=x_1y_7+...+x_7y_1.
\end{equation}

Denote by $\Sym^2(k^7)$ the space of such self-adjoint matrices. Consider the particular element
\[
A:=\matr{
	0&0&1&0&0&1&0\\
	1&0&0&0&0&0&1\\
	0&1&0&0&0&0&0\\
	0&0&0&0&0&0&0\\
	0&1&0&0&0&0&1\\
	0&0&1&0&1&0&0\\
	0&0&0&0&0&1&0
}\in\Sym^2(k^7).
\]
It is regular semisimple, meaning that $A$ has $7$ distinct eigenvalues in some finite extension of $k$. Fix any non-trivial additive character $\psi:k\ra\Cc$. The map $B\mapsto \psi(\Tr(AB))$ defines a character on the abelian group $\Sym^2(k^7)$. We denote by $\phi_{A}$ the composition $\phi_{A}:G(F)_{\frac{1}{2}}\sra G(F)_{\frac{1}{2}}/G(F)_1\cong\Sym^2(k^7)\xra{B\mapsto \psi(\Tr(AB))}\Cc$. Our main result is
\begin{theorem}\label{main}
	Any irreducible component of the compact induction
	\[
	\operatorname{c-ind}_{G(F)_{\frac{1}{2}}}^{G(F)}\phi_{A}
	\]
	is a supercuspidal representation of $G(F)$ whose wave-front set contains an orbit of Jordan type $(43)$ and another orbit of Jordan type $(511)$. In particular, its geometric wave-front set is not a singleton.
\end{theorem}

The supercuspidal representations are the so-called epipelagic representations in \cite{RY14}. Similar compact induction from $G(F)_{n+\frac{1}{2}}$ for any $n\in\Z_{\ge 0}$ should have the same wave-front set, but the proof requires longer harmonic analysis so we limit our discussion to $n=0$. The local character expansions of these representations were studied in \cite[\S6]{Ts17}; many ideas are picked up from there. 

The paper will be structured as follows: we explain the reasoning behind the choice of $A$ in Section \ref{geom}, and prove Theorem \ref{main} in Section \ref{harm}. Additionally, in Section \ref{mot} we discuss the conceptual explanation and motivation for our construction.

\subsection*{Acknowledgment} The author wishes to express gratitude to his advisor, Benedict Gross, who suggested studying the local character expansion of these epipelagic representations as an initial thesis problem back in 2012, and gave him invaluable guidance. Around the same time, Xiaoheng Jerry Wang taught the author a lot about rationality questions, which continues to influence his approach to these topics including in this paper. For this the author sincerely appreciates. The author is grateful to Sug Woo Shin for stimulating discussions, and to Jean-Loup Waldspurger for insightful comments and for pointing out a significant mistake in a draft of this paper. Additionally, the author is very grateful to the organizers of the 2022 Summer School on the Langlands program and the IHES for an extremely stimulating and wonderful environment; this work was completed during the summer school. Lastly, the author wishes to thank the referee(s) and Chi-Heng Lo for very helpful comments and suggestions on previous versions, as well as ChatGPT for polishing some English writings.\p

\section{The choice of $A$}\label{geom}

The map
\[
G(F)_0/G(F)_{\frac{1}{2}}\ira GL_7(\CO_E)/(\Id_7+\mathfrak{m}_EM_{7\times 7}(E))\cong GL_7(k)
\]
identifies $G(F)_0/G(F)_{\frac{1}{2}}$ as a subgroup of $GL_7(k)$. This subgroup is the group $O_7(k)$ of orthogonal matrices on the quadratic space $k^7$ with respect to the form (\ref{form2}). The group $G(F)_0/G(F)_{\frac{1}{2}}\cong O_7(k)$ acts by conjugation on $G(F)_{n+\frac{1}{2}}/G(F)_{n+1}\cong \Sym^2(k^7)$ and this action expectedly is the natural conjugation action of orthogonal matrices on self-adjoint matrices. What we need about the matrix $A$ is outlined in the following four lemmas:

\begin{lemma}\label{stabilizer} The stabilizer subgroup scheme of $A$ in $O_7$ is abelian, $2$-torsion of order $2^7$.
\end{lemma}

\begin{proof} An operator $g$ is orthogonal and commutes with a given regular semisimple self-adjoint operator if and only if $g$ acts as $\pm1$ on each of its eigenspaces, hence the result. In fact, in \cite[\S3.1]{Ts17} it is explained that this group mod $\pm1$ is the $2$-torsion of a Jacobian of a genus $3$ hyperelliptic curve.
\end{proof}

\begin{lemma}\label{hard} Any $O_7(k)$-conjugate of $A$ is {\bf NOT} of the following shape
\[
\matr{
	*&*&*&*&*&*&*\\
	*&*&*&*&*&*&*\\
	0&*&*&*&*&*&*\\
	0&*&*&*&*&*&*\\
	0&0&*&*&*&*&*\\
	0&0&0&*&*&*&*\\
	0&0&0&0&0&*&*
}
\]
where the symbol $*$ indicates that the entry can be anything in $k$, and the symbol $0$ indicates that the entry has to be zero.
\end{lemma}

Replacing one pair of $0$ by $*$, we would like

\begin{lemma}\label{triv} There exists an $O_7(k)$-conjugate of $A$ of the following shape
	\[
	\matr{
		*&*&*&*&*&*&*\\
		*&*&*&*&*&*&*\\
		0&*&*&*&*&*&*\\
		0&*&*&*&*&*&*\\
		0&*&*&*&*&*&*\\
		0&0&*&*&*&*&*\\
		0&0&0&0&0&*&*
	}\text{, or even of the shape}
\matr{
	*&*&*&*&*&*&*\\
	k^{\times}&*&*&*&*&*&*\\
	0&k^{\times}&*&*&*&*&*\\
	0&0&*&*&*&*&*\\
	0&k^{\times}&*&*&*&*&*\\
	0&0&k^{\times}&0&k^{\times}&*&*\\
	0&0&0&0&0&k^{\times}&*
}
\]
where the symbol $k^{\times}$ indicates that any non-zero entry is allowed.
\end{lemma}

\begin{lemma}\label{Jac} There exists an $O_7(k)$-conjugate of $A$ of the following shape
	\[
	\matr{
		*&*&*&*&*&*&*\\
		*&*&*&*&*&*&*\\
		*&*&*&*&*&*&*\\
		0&*&*&*&*&*&*\\
		0&0&*&*&*&*&*\\
		0&0&0&*&*&*&*\\
		0&0&0&0&*&*&*
	}\text{, or even of the shape}
\matr{
*&*&*&*&*&*&*\\
*&*&*&*&*&*&*\\
k^{\times}&*&*&*&*&*&*\\
0&k^{\times}&*&*&*&*&*\\
0&0&k^{\times}&*&*&*&*\\
0&0&0&k^{\times}&*&*&*\\
0&0&0&0&k^{\times}&*&*
}
	\]
\end{lemma}

Ultimately, Lemma \ref{hard} accounts for the exclusion of nilpotent orbits of Jordan type $(52)$ and anything larger from the wave-front set, while Lemma \ref{triv} (respectively, Lemma \ref{Jac}) ensures the inclusion of a nilpotent orbit of Jordan type $(511)$ (respectively, type $(43)$) within the wave-front set. Lemma \ref{triv} is made straightforward, as we choose $A$ to have the required shape. The other two lemmas can be verified using computer checks, as explained below:

\begin{proof}[Proof of Lemma \ref{hard}]
	Since $O_7=SO_7\times\{\pm\mathrm{Id}_7\}$ and $-\mathrm{Id}_7$ acts trivially, there is no harm in replacing $O_7$ by $SO_7$.
	Let $B(k)\subset SO_7(k)$ be the Borel subgroup (thanks to (\ref{form2})) of upper triangular orthogonal matrices. Since the shape we need to exclude is preserved under conjugation by $B(k)$, Lemma \ref{hard} is equivalent to the emptyness of the following set:
	\begin{equation}\label{Sym2cover}
	X^{(\ref{hard})}_A(k):=\{g\in SO_7(k)/B(k)\;|\;g^{-1}Ag\in
	\matr{
		*&*&*&*&*&*&*\\
		*&*&*&*&*&*&*\\
		0&*&*&*&*&*&*\\
		0&*&*&*&*&*&*\\
		0&0&*&*&*&*&*\\
		0&0&0&*&*&*&*\\
		0&0&0&0&0&*&*
	}
	\}.
	\end{equation}
	Note that $SO_7(k)/B(k)=(SO_7/B)(k)$ is the set of isotropic flags in $k^7$, or equivalently $k$-points of the $9$-dimensional flag variety. This set is small enough to be listed by computer programs. We refer to Appendix \ref{AppReadme} for a simplified program in {\bf Magma} language, and for more examples when $p=5$ and $p=7$. We also refer to \S\ref{subsec:AG} for some conceptual explanations about this lemma.
\end{proof}

\begin{proof}[Proof of Lemma \ref{Jac}]
	We have the following choice of $g\in O_7(k)$ and the resulting conjugate:
	\[
	g=
	\matr{
	0&0&1&1&1&2&2\\
	0&1&0&1&0&1&2\\
	1&0&0&1&0&0&1\\
	1&1&0&0&2&0&2\\
	1&2&0&2&1&0&1\\
	0&1&2&2&0&1&1\\
	0&1&0&1&1&1&2
	}
	,\;g^{-1}Ag=
	\matr{
	0&0&0&0&2&1&2\\
	0&0&1&0&0&0&1\\
	1&0&1&0&0&0&2\\
	0&1&0&1&0&0&0\\
	0&0&1&0&1&1&0\\
	0&0&0&1&0&0&0\\
	0&0&0&0&1&0&0
	}.
	\]
	We include in Appendix \ref{App2.4} a program to verify these matrices.
\end{proof}

%In fact, using the method of \cite[Thm. 3.1 and Lemma 4.6]{Ts17} it is possible to show that Lemma \ref{Jac} is always valid because such conjugates correspond to rational points on a torsor of the Jacobian of the genus $3$ hyperelliptic curve over $k$ studied {\it op. cit.}, and such a torsor always have rational points due to Lang's theorem. On the other hand, Lemma \ref{hard} and \ref{triv} are concerned with rational points on some surfaces over $k$ that we haven't fully understood; we just use a computer program to find one with/without rational points.
Let us conclude with two linear algebra lemmas regarding regular self-adjoint matrices for later use.

\begin{lemma}\label{rs} A regular semisimple self-adjoint matrix in $\Sym^2(k^7)$ cannot be (properly) blockwise triangular.
\end{lemma}

\begin{proof} Suppose on the contrary that some regular semisimple self-adjoint matrix is, for example, of the form
\[
\matr{
a&b&*&*&*&*&*\\
c&d&*&*&*&*&*\\
0&0&*&*&*&*&*\\
0&0&*&*&*&*&*\\
0&0&*&*&*&*&*\\
0&0&0&0&0&d&b\\
0&0&0&0&0&c&a\\
}
\]
Then the top-left $\matr{a&b\\c&d}$ and the bottom-right $\matr{d&b\\c&a}$ have the same eigenvalues and thus the matrix is not regular semisimple. The same reasoning works for all blockwise triangular shapes.
\end{proof}

\begin{lemma}\label{conj} Any $O_7(k)$-conjugate of $A$ is not of the following shape
\[
\matr{
	*&*&*&*&*&*&*\\
	*&*&*&*&*&*&*\\
	*&*&*&*&*&*&*\\
	*&*&*&*&*&*&*\\
	0&0&0&*&*&*&*\\
	0&0&0&*&*&*&*\\
	0&0&0&*&*&*&*
}
\]
\end{lemma}

\begin{proof} We claim that any $M\in\Sym^2(k^7)$ of the shape above can be conjugate into the shape below
	\[
	\matr{
		*&*&*&*&*&*&*\\
		*&*&*&*&*&*&*\\
		0&*&*&*&*&*&*\\
		0&0&*&*&*&*&*\\
		0&0&0&*&*&*&*\\
		0&0&0&0&*&*&*\\
		0&0&0&0&0&*&*
	}
	\]
by an element $\til{C}:=\matr{C&&\\&1&\\&&(C^t)^{-1}}\in O_7(k)$ where $C$ is an invertible $3\times 3$ matrix and $C^t$ is the transpose with respect to the anti-diagonal. The last displayed shape is then forbidden by Lemma \ref{hard}. To prove the claim, denote by $r_i(M)$ the first $3$ entries of the $i$-th row of $M$, as a $3$-dimensional row vector. Conjugation (from the right) by $\til{C}$ changes $r_i(M)$ to $r_i(M)C$ (and symmetrically so for the last $3$ entries of the last four columns) and preserves the lower-bottom $3\times3$ of 0's. The claim is thus equivalent to that $r_4(M)C$ is of the shape $[0\;0\;*]$ and $r_3(M)C$ is of the shape $[0\;*\;*]$. This is always achieved by some $C$ for any two row vectors.
\end{proof}

\section{Harmonic analysis}\label{harm}

Firstly let us specify our Fourier transform. We choose an additive character $\psi:F\ra\Cc$ with kernel equal to the maximal ideal $\mathfrak{m}_F=3\Z_3$, so that $\psi$ induces the identically named character on $k$ used in the introduction. The Lie algebra $\Lg(F)$ is the space of anti-hermitian $7\times 7$ matrices with respect to (\ref{form1}), i.e.
\[
\Lg(F)=\{T\in M_{7\times 7}(E)\;|\;\langle Tx,y\rangle=\langle x,-Ty\rangle,\;\forall x,y\in E^7\}.
\]
We have a pairing $\beta:\Lg(F)\times\Lg(F)\ra F$ given by $\beta(X,Y)=\Tr(XY)$, with which we define $\hat{f}$ for $f\in C_c^{\infty}(\Lg(F))$ as
\[
\hat{f}(X):=\int_{\Lg(F)}\psi(\beta(X,Y))f(Y)dY.
\]
Any positive translation-invariant measure on $\Lg(F)$ suffices. We have $\hat{\hat{f}}=C.f^{\#}$ for some constant $C\in\R_{>0}$, where $f^{\#}(X)=f(-X)$. The number $C$ depends on chosen measure. Since we ultimately care only about whether some coefficients are non-zero, and they all scale with $C$, the exact choice won't matter. This implies that for any distribution $D\in C_c^{\infty}(\Lg(F))^*$, we have $\hat{D}(\hat{f})=D(\hat{\hat{f}})$ (by definition) is equal to $C\cdot D(f^{\#})$. When $f=f^{\#}$, which will always be the case below, we have $\hat{D}(\hat{f})=C\cdot D(f)$.

The Lie algebra $\Lg(F)$ has a decreasing Moy-Prasad filtration $(\Lg(F)_r)_{r\in\frac{1}{2}\Z}$ where $\Lg(F)_r$ consists of matrices in $\Lg(F)$ whose entries live in $\mathfrak{m}_E^{2r}$. They satisfy $\Lg(F)_{r+n}=p^n\Lg(F)_r$ for any $n\in\Z$. Moreover, for any $n\in\Z$, dividing by $\varpi^{2n+1}$ gives a map
\begin{equation}\label{map}
	\Lg(F)_{n+\frac{1}{2}}/\Lg(F)_{n+1}\ira \mathfrak{m}_E^{2n+1}M_{7\times 7}(\CO_E)/\mathfrak{m}_E^{2n+2}M_{7\times 7}(\CO_E)\cong M_{7\times 7}(k).
\end{equation}
The map again identifies $\Lg(F)_{n+\frac{1}{2}}/\Lg(F)_{n+1}$ as the space $\Sym^2(k^7)$. We note that for any $r\in\frac{1}{2}\Z$ with $r\ge 1$, we have that the exponential map $\exp$ defines an isomorphism between $\Lg(F)_r$ and $G(F)_r$ that transports (\ref{map}) to (\ref{map1}). In fact, we have the Cayley transform $\se(X):=(\Id+\frac{1}{2}X)(\Id-\frac{1}{2}X)^{-1}$ that gives an isomorphism $\Lg(F)_r\xra{\sim}G(F)_r$ for any $r\ge\frac{1}{2}$, with inverse $\se^{-1}(g):=2(g-\Id)(g+\Id)^{-1}$. It is a general fact that in the local character expansion (\ref{LCE}) one can replace $\exp$ by $\se$ to get the same expansion with same coefficients. For our purpose, we will only plug in specific test functions (to be introduced right below) on $\Lg(F)$ whose composition with $\log$ will be evidently the same as that with $\se^{-1}$. Hence we can and will replace all potential $\exp$ by $\se$ and $\log$ by $\se^{-1}$.

For any subset $S\subset\Sym^2(k^7)$, let us denote by $f_n^S$ the function supported on $\Lg(F)_{-n-\frac{1}{2}}$ whose value is $1$ at those elements whose image under $\Lg(F)_{-n-\frac{1}{2}}\sra\Lg(F)_{-n-\frac{1}{2}}/\Lg(F)_{-n}\xra[\sim]{(\ref{map})}\Sym^2(k^7)$ goes to $S$, and $0$ otherwise. We remark that the support of their Fourier transform $\hat{f}_n^S$ will be contained in $\Lg(F)_{n+\frac{1}{2}}$. Consider subsets $S_{(7)},S_{(61)},S_{(52)},S_{(43)},S_{(511)}\subset\Sym^2(k^7)$ as
\[
S_{(7)}=\{\matr{
	0&0&0&0&0&0&0\\
	k^{\times}&0&0&0&0&0&0\\
	0&k^{\times}&0&0&0&0&0\\
	0&0&k^{\times}&0&0&0&0\\
	0&0&0&k^{\times}&0&0&0\\
	0&0&0&0&k^{\times}&0&0\\
	0&0&0&0&0&k^{\times}&0
}\},\;
S_{(61)}=\{\matr{
	0&0&0&0&0&0&0\\
	k^{\times}&0&0&0&0&0&0\\
	0&k^{\times}&0&0&0&0&0\\
	0&0&0&0&0&0&0\\
	0&0&k^{\times}&0&0&0&0\\
	0&0&0&0&k^{\times}&0&0\\
	0&0&0&0&0&k^{\times}&0
}\},
\]
\[
S_{(52)}=\{\matr{
	0&0&0&0&0&0&0\\
	k^{\times}&0&0&0&0&0&0\\
	0&0&0&0&0&0&0\\
	0&k^{\times}&0&0&0&0&0\\
	0&0&k^{\times}&0&0&0&0\\
	0&0&0&k^{\times}&0&0&0\\
	0&0&0&0&0&k^{\times}&0
}\},\;
S_{(43)}=\{\matr{
	0&0&0&0&0&0&0\\
	0&0&0&0&0&0&0\\
	k^{\times}&0&0&0&0&0&0\\
	0&k^{\times}&0&0&0&0&0\\
	0&0&k^{\times}&0&0&0&0\\
	0&0&0&k^{\times}&0&0&0\\
	0&0&0&0&k^{\times}&0&0
}\},
\]
\[
S_{(511)}=\{\matr{
	0&0&0&0&0&0&0\\
	k^{\times}&0&0&0&0&0&0\\
	0&k^{\times}&0&0&0&0&0\\
	0&0&0&0&0&0&0\\
	0&k^{\times}&0&0&0&0&0\\
	0&0&k^{\times}&0&k^{\times}&0&0\\
	0&0&0&0&0&k^{\times}&0
}\},
\]
where the symbol $k^{\times}$ indicates that the entry can be any element in $k^{\times}=\F_3^{\times}$. In particular, one sees that an element in each $S_{\lambda}$ is nilpotent with Jordan type $\lambda$. Theorem \ref{main} will be proved by using the Fourier transforms $\hat{f}_n^S$ as the test functions. The required calculation is given in the following two propositions 

\begin{proposition}\label{ss} Let $\pi\subset\operatorname{c-ind}_{G(F)_{\frac{1}{2}}}^{G(F)}\phi_A$ be any irreducible component. For any integer $n\in\Z_{>0}$ we have
\[
(\Theta_{\pi}\circ\se)(\hat{f}_{n}^{S_{(7)}})=(\Theta_{\pi}\circ\se)(\hat{f}_{n}^{S_{(61)}})=(\Theta_{\pi}\circ\se)(\hat{f}_{n}^{S_{(52)}})=0
\]
and
\[
(\Theta_{\pi}\circ\se)(\hat{f}_{n}^{S_{(43)}})>0,\;(\Theta_{\pi}\circ\se)(\hat{f}_{n}^{S_{(511)}})>0.
\]
where $(\Theta_{\pi}\circ\se)$ denotes the pullback of $\Theta_{\pi}|_{G(F)_{n+\frac{1}{2}}}$ to $\Lg(F)_{n+\frac{1}{2}}$.
\end{proposition}

We postpone the longer proof of the above proposition to the end of the section. Meanwhile, our group $U_7(E/F)$ has a unique nilpotent orbit $\CO_7$ of Jordan type $(7)$, two nilpotent orbits $\CO_{61,+}$ and $\CO_{61,-}$ of Jordan type $(61)$, two nilpotent orbits $\CO_{52,+}$ and $\CO_{52,-}$ of Jordan type $(52)$, two nilpotent orbits $\CO_{43,+}$ and $\CO_{43,-}$ of Jordan type $(43)$, and two nilpotent orbits $\CO_{511,+}$ and $\CO_{511,-}$ of Jordan type $(511)$. See e.g. \cite[\S4]{Ts17}. We have

\begin{proposition}\label{nil} For any $n\in\Z$, we have
\[
\{\CO\text{ nilpotent orbit }|\;I_{\CO}(f^{S_{(7)}}_n)\not=0\}=\{\CO_7\}.
\]
\[
\{\CO\text{ nilpotent orbit }|\;I_{\CO}(f^{S_{(61)}}_n)\not=0\}=\{\CO_7,\CO_{61,+},\CO_{61,-}\}.
\]
\[
\{\CO\text{ nilpotent orbit }|\;I_{\CO}(f^{S_{(52)}}_n)\not=0\}=\{\CO_7,\CO_{61,+},\CO_{61,-},\CO_{52,+},\CO_{52,-}\}.
\]
\[
\{\CO\text{ nilpotent orbit }|\;I_{\CO}(f^{S_{(43)}}_n)\not=0\}=\{\CO_7,\CO_{61,+},\CO_{61,-},\CO_{52,+},\CO_{52,-},\CO_{43,+},\CO_{43,-}\}.
\]
\[
\{\CO\text{ nilpotent orbit }|\;I_{\CO}(f^{S_{(511)}}_n)\not=0\}=\{\CO_7,\CO_{61,+},\CO_{61,-},\CO_{52,+},\CO_{52,-},\CO_{511,+},\CO_{511,-}\}.
\]
In other words, for each partition $\lambda$ appearing above, those nilpotent orbits in $\Lg(F)$ such that $I_{\CO}(f_n^{S_{\lambda}})\not=0$ are exactly those with Jordan type $\lambda$ or larger. We remark that whenever $I_{\CO}(f_n^S)\not=0$, it is by definition positive.
\end{proposition}

\begin{proof}[Proof of Proposition \ref{nil}] We first prove the $\subset$ direction. To say $\CO$ is in the set on the LHS is to say that some element $e\in\CO$ lives in the support of the function $f_n^S$, which implies that as a $7\times 7$ matrix over $E$, the reduction $\bar{e}$ of $e$ mod $\mathfrak{m}_E^{-2n}$ is nilpotent of the Jordan type indicated. Since $e$ itself is also nilpotent, the Jordan type for $e$ can only be larger and hence $\CO\ni e$ belongs to the RHS. Indeed, a nilpotent matrix $\bar{e}$ has Jordan type $\lambda$ or larger iff for the dual partition $\lambda^t=(\ell_1\ge\ell_2\ge ...\ell_s)$ we have $\operatorname{nullity}(\bar{e}^i)\le \ell_1+...+\ell_i$. The assertion direction then follows from that $\operatorname{nullity}(e^i)\le \operatorname{nullity}(\bar{e}^i)$.
	
For the $\supset$ direction, one has to find $e\in\CO\cap\supp(f_n^S)$ for each orbit $\CO$ on the RHS. Since $\supp(f_n^S)$ is by definition closed, it suffices to find such $e$ for each $\CO$ minimal on the RHS. Thus the proposition is proved by finding $e_{7}\in\CO_{7}\cap\supp(f_n^{S_{(7)}})$ and  $e_{\lambda,\pm}\in\CO_{\lambda,\pm}\cap\supp(f_n^{S_\lambda})$ for partitions $\lambda\in\{(61),(52),(43),(511)\}$. Let $c:=\varpi^{-2n-1}$, $d_+\in\CO_F^{\times}$ be any square and $d_-\in\CO_F^{\times}$ be any non-square. The desired nilpotent elements can be given by
\[{\tiny\hskip-.5cm
e_{7}=c\matr{
	0&0&0&0&0&0&0\\
	1&0&0&0&0&0&0\\
	0&1&0&0&0&0&0\\
	0&0&1&0&0&0&0\\
	0&0&0&1&0&0&0\\
	0&0&0&0&1&0&0\\
	0&0&0&0&0&1&0
},\;e_{61,\pm}=c\matr{
0&0&0&0&0&0&0\\
1&0&0&0&0&0&0\\
0&1&0&0&0&0&0\\
0&0&0&0&0&0&0\\
0&0&d_{\pm}&0&0&0&0\\
0&0&0&0&1&0&0\\
0&0&0&0&0&1&0
},\;e_{52,\pm}=c\matr{
0&0&0&0&0&0&0\\
1&0&0&0&0&0&0\\
0&0&0&0&0&0&0\\
0&1&0&0&0&0&0\\
0&0&d_{\pm}&0&0&0&0\\
0&0&0&1&0&0&0\\
0&0&0&0&0&1&0
},}
\]
\[{\tiny
e_{43,\pm}=c\matr{
	0&0&0&0&0&0&0\\
	0&0&0&0&0&0&0\\
	1&0&0&0&0&0&0\\
	0&1&0&0&0&0&0\\
	0&0&d_{\pm}&0&0&0&0\\
	0&0&0&1&0&0&0\\
	0&0&0&0&1&0&0
},\;e_{511,\pm}=c\matr{
0&0&0&0&0&0&0\\
1&0&0&0&0&0&0\\
0&1/2&0&0&0&0&0\\
0&0&0&0&0&0&0\\
0&d_{\pm}&0&0&0&0&0\\
0&0&d_{\pm}&0&1/2&0&0\\
0&0&0&0&0&1&0
}}.
\]
\vskip-.5cm
\end{proof}

\begin{proof}[Proof of Theorem \ref{main}] Choose an $n\in\Z_{>0}$ large enough so that (\ref{LCE}) is valid on $\se(\Lg(F)_{n+\frac{1}{2}})=G(F)_{n+\frac{1}{2}}$. We plug $\hat{f}_n^{S_{\lambda}}\in C^{\infty}(\Lg(F)_{n+\frac{1}{2}})$ via $\se$ for $\lambda\in\{(7),(61),(52),(43),(511)\}$ into (\ref{LCE}), namely
\begin{equation}\label{apply}
	(\Theta_{\pi}\circ\se)(\hat{f}_n^{S_{\lambda}})=\sum_{\CO}c_{\CO}(\pi)\cdot\hat{I}_{\CO}(\hat{f}_n^{S_{\lambda}})=\sum_{\CO}c_{\CO}(\pi)\cdot I_{\CO}(f_n^{S_{\lambda}})
\end{equation}
That $(\Theta_{\pi}\circ\se)(\hat{f}_{n}^{S_{(7)}})=0$ and the first line of Proposition \ref{nil} asserts that $c_{\CO_7}(\pi)\cdot I_{\CO_7}(f_n^{S_{(7)}})=0$ and $I_{\CO_7}(f_n^{S_{(7)}})>0$, and thus $c_{\CO_7}(\pi)=0$. Next we look at $(\Theta_{\pi}\circ\se)(\hat{f}_{n}^{S_{(61)}})=0$. Combining with the second line of Proposition \ref{nil} we have $c_{\CO_{61,+}}(\pi)\cdot I_{\CO_{61,+}}(f_n^{S_{(61)}})+c_{\CO_{61,-}}(\pi)\cdot I_{\CO_{61,-}}(f_n^{S_{(61)}})=0$. By \cite[Cor. 1.17]{MW87}, that $c_{\CO_7}=0$ implies that the numbers $c_{\CO_{61,\pm}}\ge 0$, because up to a normalizing positive constant it is the dimension of specific degenerate Whittaker models. Since $I_{\CO_{61,\pm}}(f_n^{S_{(61)}})>0$. This implies that $c_{\CO_{61,+}}(\pi)=c_{\CO_{61,-}}(\pi)=0$. Continue with $(\Theta_{\pi}\circ\se)(\hat{f}_{n}^{S_{(52)}})=0$ we get $c_{\CO_{52,+}}(\pi)\cdot I_{\CO_{52,+}}(f_n^{S_{(52)}})+c_{\CO_{52,-}}(\pi)\cdot I_{\CO_{52,-}}(f_n^{S_{(52)}})=0$ and consequently $c_{\CO_{52,+}}(\pi)=c_{\CO_{52,-}}(\pi)=0$.
 
Next, we look at $(\Theta_{\pi}\circ\se)(\hat{f}_{n}^{S_{(43)}})>0$. Thanks to the fourth line of Proposition \ref{nil}, from the RHS of (\ref{apply}) we have $c_{\CO_{43,+}}(\pi)\cdot I_{\CO_{43,+}}(f_n^{S_{(43)}})+c_{\CO_{43,-}}(\pi)\cdot I_{\CO_{43,-}}(f_n^{S_{(43)}})>0$, and thus at least one of $c_{\CO_{43,\pm}}(\pi)\not=0$ (actually positive by \cite[Cor. 1.17]{MW87}). Likewise, $(\Theta_{\pi}\circ\se)(\hat{f}_{n}^{S_{(511)}})>0$ and the fifth line of Proposition \ref{nil} gives that one of $c_{\CO_{511,\pm}}(\pi)\not=0$. This proves our main theorem except for the statement that the components are supercuspidal, which is part of \cite[Prop. 2.4]{RY14}.

(We thank an anonymous referee for pointing out the need and use of \cite{MW87}.)
\end{proof}

\begin{proof}[Proof of Proposition \ref{ss}] In the rest of this section, whenever $X$ is some object on which $G(F)$ acts on the left (typically by conjugation) we will denote by ${}^gX$ the left action and $X^g:={}^{g^{-1}}X$.

Let us first recall the structure of $\operatorname{c-ind}_{G(F)_{\frac{1}{2}}}^{G(F)}\phi_A$. By \cite[Prop. 2.4]{RY14},  $\operatorname{c-ind}_{G(F)_0}^{G(F)}\til{\phi}_A$ is a direct sum of finitely many irreducible supercuspidal representations. Each of them is of the form $\operatorname{c-ind}_{G(F)_0}^{G(F)}\til{\phi}_A$ for some $\til{\phi}_A\in\Irr(G(F)_0)$ %that is $\phi_A$-isotypic\footnote{Here we abuse the language and say that a representation $\rho$ of $G(F)_0$ is $\eta$-isotypic for $\eta\in\Irr(G(F)_{\frac{1}{2}})$ if 
such that $\til{\phi}_A|_{G(F)_{\frac{1}{2}}}$ is a direct sum of a finite number of $G(F)_0$-conjugates of $\phi_A$. For any $f\in C_c^{\infty}(G(F))$, definition of compact induction gives the character as
\begin{equation}\label{HC}
\Theta_{\pi}(f)=\sum_{g\in G(F)_0\bsl G(F)}\langle\Theta_{\til{\phi}_A},{}^gf\rangle.
\end{equation}
Moreover, in \cite[Prop. 2.4]{RY14}, the representation $\til{\phi}_A$ is induced from $\operatorname{Stab}_{G(F)_0}(A)$, the preimage in $G(F)_0$ of the stabilizer of $A$ in $G(F)_0/G(F)_{\frac{1}{2}}=O_7(k)$. In particular the character $\Theta_{{}^g\til{\phi}_A}$ is supported on $G(F)_0$-conjugates of $\operatorname{Stab}_{G(F)_0}(A)$. By Lemma \ref{stabilizer}, any element in $\operatorname{Stab}_{G(F)_0}(A)$ is either in $G(F)_{\frac{1}{2}}$, or has an eigenvalue $\lambda$ with $\operatorname{val}(\lambda-1)=0$ (in fact $|\lambda+1|<1$). In particular, the only elements in $\operatorname{Stab}_{G(F)_0}(A)$ that can meet conjugates of elements in $G(F)_{\frac{1}{2}}$ are those in $G(F)_{\frac{1}{2}}$ themselves. That is to say, if $f$ is supported on $G(F)_{\frac{1}{2}}$, then (\ref{HC}) is simplified to
\begin{equation}\label{HC2}
	\Theta_{\pi}(f)=\sum_{g\in G(F)_0\bsl G(F)}\langle\Theta_{\til{\phi}_A}|_{G(F)_{\frac{1}{2}}},{}^gf\rangle.
\end{equation}
Since $\til{\phi}_A|_{G(F)_{\frac{1}{2}}}$ is a direct sum of $G(F)_0$-conjugates of $\phi_A$, the restriction $\Theta_{\til{\phi}_A}|_{G(F)_{\frac{1}{2}}}$ is a multiple of  $\Theta_{\phi_A}$ average by $G(F)_0/G(F)_{\frac{1}{2}}$. More precisely, let $\bar{f}_A$ be the function on $\Sym^2(k^7)$ defined by
\[
\bar{f}_A(B)=\#\{\bar{g}\in O_7(k)\;|\;\Ad(\bar{g})(B)=A\}
\]
and let $f_A$ denote the pullback of $\bar{f}_A$ under the map $\Lg(F)_{-\frac{1}{2}}/\Lg(F)_0\sra\Sym^2(k^7)$. By construction, $f_A$ is invariant under conjugation by $G(F)_0$. We then have
\begin{equation}\label{char}
(\Theta_{\til{\phi}_A}|_{G(F)_{\frac{1}{2}}})\circ\se=C\cdot\hat{f}_A
\end{equation}
for some constant $C>0$.

We have the Cartan decomposition that
	\begin{equation}\label{Cartan}
		G(F)_0\bsl G(F)=\bigsqcup_{d\in D}G(F)_0\bsl G(F)_0\cdot d\cdot G(F)_0
	\end{equation}
where
\[
\hskip-.25cm
D=\{
\matr{
	\varpi^{d_3}&&&&&&\\
	&\varpi^{d_2}&&&&&\\
	&&\varpi^{d_1}&&&&\\
	&&&1&&&\\
	&&&&(-\varpi)^{-d_1}&&\\
	&&&&&(-\varpi)^{-d_2}&\\
	&&&&&&(\varpi)^{-d_3}
}\;|\;d_3\ge d_2\ge d_1\ge 0\text{ are integers}.\}
\]

Combining with (\ref{HC2}) and (\ref{char}), this gives for $f'\in C^{\infty}(\Lg(F)_{\frac{1}{2}})$ and $f=f'\circ\se^{-1}$ that
\[
\Theta_{\pi}(f)=\sum_{d\in D}\left(\sum_{g\in (G(F)_0^d\cap G(F)_0)\bsl G(F)_0}\langle C.\hat{f}_A,{}^{dg} f'\rangle\right).
\]
\begin{equation}\label{HC3}
=C\cdot\sum_{d\in D}\left(\sum_{g\in (G(F)_0^d\cap G(F)_0)\bsl G(F)_0}\langle (f_A)^{d},{}^g\hat{f'}\rangle\right)
\end{equation}

Suppose $f'=\hat{f}_n^{S_{(43)}}$ so that $\hat{f'}=f_n^{S_{(43)}}$ up to some positive constant. Take $g=\Id$, $d_3=3n$, $d_2=2n$ and $d_1=n$. The intersection $\supp(f_A)^{d}\cap \supp(f_n^{S_{(43)}})\not=\emptyset$ thanks to Lemma \ref{Jac}. Hence $(\Theta_{\pi}\circ\se)(\hat{f}_n^{S_{(43)}})>0$. Likewise when $g=\Id$, $d_3=4n$, $d_2=2n$ and $d_1=0$, thanks to Lemma \ref{triv} we have $\supp(f_A)^{d}\cap \supp(f_n^{S_{(511)}})\not=\emptyset$ and thus $(\Theta_{\pi}\circ\se)(\hat{f}_n^{S_{(511)}})>0$. 

Suppose $f'=\hat{f}_n^{S_{\lambda}}$ for some $\lambda\in\{(7),(61),(52)\}$ so that $\hat{f'}=f_n^{S_{\lambda}}$ up to some positive constant. We want to prove that (\ref{HC3}) is zero. (We will see that $\lambda=(52)$ is the essential case.) In the sum (\ref{HC3}), we are only concerned with those $(d,g)$ for which the conjugate $\supp((f_A)^d)=\supp(f_A)^d$ meets $\supp({}^gf_n^{S_\lambda})\subset\Lg(F)_{-n-\frac{1}{2}}$. In other words, what could contribute is $\til{B}\in\supp(f_A)\subset\Lg(F)_{-\frac{1}{2}}$ such that $\til{B}^d\in\supp({}^gf_n^{S_\lambda})\subset\Lg(F)_{-n-\frac{1}{2}}$. We will show that such $d\in D$ and $\til{B}$ don't exist.

For a matrix $\til{B}=(B_{ij})_{-3\le i,j\le 3}\in\Lg(F)\subset M_{7\times 7}(E)$ (with the special indices), we have $(\til{B}^d)_{ij}=\pm\varpi^{d_i-d_j}B_{ij}$ where we write $d_{-i}=-d_i$ for $i=1,2,3$ and $d_0=0$. Suppose $d_3-d_2>2n$. In this case the conditions $(\til{B}^d)_{3,j}\in\mathfrak{m}_E^{-2n-1}$ for $-3\le j\le 2$ implies that the reduction $B\in \Lg(F)_{-\frac{1}{2}}/\Lg(F)_{0}\cong\Sym^2(k^7)$ of $\til{B}$ has to be of the form
\[
\matr{
*&*&*&*&*&*&*\\
0&*&*&*&*&*&*\\
0&*&*&*&*&*&*\\
0&*&*&*&*&*&*\\
0&*&*&*&*&*&*\\
0&*&*&*&*&*&*\\
0&0&0&0&0&0&*
}
\]
By construction of $f_A$, we need $B$ to be an $O_7(k)$-conjugate of $A$. Thus the above shape is not possible by Lemma \ref{rs}. Hence $d_3-d_2\le 2n$. Similarly Lemma \ref{rs} gives $d_2-d_1\le 2n$ and $d_1\le 2n$. In fact we have $d_1\le n$, for otherwise $2d_1>2n$ and $B$ has to be of the shape in Lemma \ref{conj}. Next we claim that $d_3-d_2<2n$ is also not possible. Suppose $d_3-d_2<2n$. Then $(\til{B}^d)_{32}\in\varpi^{1-2n}\mathfrak{m}_E^{-1}=\mathfrak{m}_E^{-2n}$. For this and analogous reasonings for other entries, the reduction of $\til{B}^d$ in $\Lg(F)_{-n-\frac{1}{2}}/\Lg(F)_{-n}\cong\Sym^2(k^7)$ is of the following shape
\[
\matr{
	0&0&0&0&0&0&0\\
	\mathbf{0}&0&0&0&0&0&0\\
	*&*&0&0&0&0&0\\
	*&*&\mathbf{a}&0&0&0&0\\
	*&*&*&\mathbf{a}&0&0&0\\
	*&*&*&*&*&0&0\\
	*&*&*&*&*&\mathbf{0}&0
}
\]
Such a shape cannot have Jordan type (7) or (61) because of the vanishing of the highlighted spot, thus not possible for $S_{(7)}$ nor $S_{(61)}$. For the case of $S_{(52)}$, the above shape has Jordan type (52) only when $\mathbf{a}\in k^{\times}$ is non-zero. But this is only possible if $d_1\ge 2n$ which contradicts with $d_1\le n$.

We have proved that $d_3-d_2=2n$ and $d_1\le n$ are necessary conditions to yield a non-zero contribution to (\ref{HC3}). When $2d_1<2n$, the reduction of $\til{B}^d$ has the form
\[
\matr{
	0&0&0&0&0&0&0\\
	*&0&0&0&0&0&0\\
	*&*&0&0&0&0&0\\
	*&*&0&0&0&0&0\\
	*&*&\mathbf{0}&0&0&0&0\\
	*&*&*&*&*&0&0\\
	*&*&*&*&*&*&0\\
}
\]
which contains no nilpotent elements of Jordan type (52) or larger. Therefore, the only possibility is $d_1=n$. In the case $d_2=d_1=n$, the reduction of $\til{B}^d$ takes the following form
\[
\matr{
	0&0&0&0&0&0&0\\
	*&0&0&0&0&0&0\\
	*&\mathbf{0}&0&0&0&0&0\\
	*&\mathbf{0}&0&0&0&0&0\\
	*&*&*&0&0&0&0\\
	*&*&*&\mathbf{0}&\mathbf{0}&0&0\\
	*&*&*&*&*&*&0\\
}
\]
which has no nilpotent element of a Jordan block of size $5$ or larger. In the remaining case $d_2>d_1=n$, we have both $d_3-d_1>2n$ and $d_2+d_1>2n$ so that the reduction $B$ of $\til{B}$ must be of the form
\[
\matr{
	*&*&*&*&*&*&*\\
	*&*&*&*&*&*&*\\
	\mathbf{0}&*&*&*&*&*&*\\
	0&*&*&*&*&*&*\\
	0&\mathbf{0}&*&*&*&*&*\\
	0&0&\mathbf{0}&*&*&*&*\\
	0&0&0&0&\mathbf{0}&*&*
}
\]
which contradicts Lemma \ref{hard}. Thus, we have considered all possibilities of $d \in D$, and none of them contributes to the sum in (\ref{HC3}) when $f'=\hat{f}_n^{S_\lambda}$ for $\lambda\in\{(7),(61),(52)\}$.
As a result, $(\Theta_{\pi}\circ\se)(\hat{f}_n^{S_{(7)}})=(\Theta_{\pi}\circ\se)(\hat{f}_n^{S_{(61)}})=(\Theta_{\pi}\circ\se)(\hat{f}_n^{S_{(52)}})=0$. This concludes the proof of Proposition \ref{ss} and, consequently, Theorem \ref{main}.
\end{proof}

\begin{remark} We may define $\CO_{511,+}$ and $\CO_{511,-}$ in the way that $\CO_{511,+}$ meets the Lie algebra of the Levi subgroup $U_5\times U_2\subset U_7$ while $\CO_{511,-}$ does not (therefore $F$-distinguished). A careful analysis of the above shows $c_{\CO_{511,+}}(\pi)=0$ while $c_{\CO_{511,-}}(\pi)>0$. Meanwhile both $c_{\CO_{43,+}}(\pi)>0$ and $c_{\CO_{43,-}}(\pi)>0$.
\end{remark}

\section{Discussion}\label{mot}

In this section, we aim to discuss the underlying philosophy and our understanding of the construction.

\subsection{Shalika germ expansion}\label{subsec:Sha} In \cite{KM03} and \cite{KM06}, it is shown that when $p\gg\rank G$, the local character of an irreducible admission representation is a linear combination of Fourier transforms of specific non-nilpotent orbital integrals. For some moderately common representations $\pi$ (e.g. regular supercuspidals in \cite{Kal19}), we obtain a single orbital integral. In other words, there exists an element $\til{A}=\til{A}_{\pi}\in\Lg$ (identifying $\Lg$ with its dual) and $C\in\R_+$, constructed from the type for $\pi$, such that $\Theta_{\pi}\equiv C\cdot \hat{I}_{\til{A}}$ on some neighborhood of the identity. For $p$ large, a variant of the Shalika germ expansion \cite[Thm. 2.1.5]{De02a} also asserts the existence of constants $s_{\CO}({\til{A}})\in\C$ for $\CO$ running over nilpotent orbits in $\Lg(F)$ satisfying $I_{\til{A}}(f)=\sum_\CO s_{\CO}({\til{A}})I_\CO(f)$ for all functions $f$ locally constant by a sufficiently large lattice. Applying Fourier transforms, we obtain
\begin{equation}\label{eq:Sha}
\hat{I}_{\til{A}}\equiv \sum s_\CO({\til{A}})\hat{I}_\CO
\end{equation}
on a sufficiently small neighborhood of $0\in\Lg(F)$. This results in $\Theta_{\pi}\equiv \sum C\cdot s_\CO({\til{A}})\hat{I}_{\CO}$ on some neighborhood, meaning these $C\cdot s_{\CO}({\til{A}})$ are precisely the coefficients $c_{\CO}(\pi)$ in (\ref{LCE}), and the wave-front set question is equivalent to the analogous wave-front set question for the so-called Shalika germs $s_{\CO}({\til{A}})$ (i.e. the set of largest $\CO$ for which $s_{\CO}({\til{A}})\not=0$). Although we do not know if all these work for small $p$ such as $p=3$, this is the starting point of our heuristic for the wave-front set. We note that in our case, $\til{A}=\til{A}_{\pi}$ is any lift in $\Lg(F)_{-\frac{1}{2}}$ of $A$ in the introduction under (\ref{map}) for $n=-1$. 

\subsection{Springer theory}\label{subsec:Springer} Determination of Shalika germs are generally very difficult. Nevertheless, an analogue for (\ref{eq:Sha}) over a finite field is well-understood through the classical Springer theory. In this case, we should consider ``a neighborhood of $0$ in a p-adic Lie algebra'' as analogous to ``the set of topologically nilpotent elements,'' and take its finite-field analogue to be ``the set of nilpotent elements.'' Suppose $H$ is a reductive group over the  residue field $k$ and $\Lh:=\Lie H$. For simplicity we assume $H$ split with $\mathrm{char}(k)\gg\rank H$. When $A\in\Lh(k)$ is regular semisimple, we have for any nilpotent $n\in\Lh^{nil}(k)$ that \cite[Thm. 4.4]{Spr76}:
\begin{equation}\label{eq:Spr}
\hat{I}_A(n)=\Tr(w.\Frob,R\pi_*\underline{\mathbb{Q}_{\ell}}_{\til{\Lh}^{nil}}|_n),
\end{equation}
where $I_A$ is the sum of values on $H(k)$-conjugates of $A$, $\hat{I}_A$ is its Fourier transform, $w$ is the Weyl group element classifying the rational conjugacy class of $Z_H(A)$, and $\pi:\til{\Lh}^{nil}=\{(A,g)\in\Lh^{nil}\times H/B\;|\;g^{-1}Ag\in\Lie B\}\ra\Lh^{nil}$ is the Springer resolution and $R\pi_*\underline{\mathbb{Q}_{\ell}}_{\til{\Lh}^{nil}}$ is the so-called Springer sheaf equipped with the Springer $W$-action. The Springer sheaf is a direct sum of various IC extensions of equivariant local systems on nilpotent orbits. The RHS of (\ref{eq:Spr}) can be rewritten in terms of class functions on nilpotent elements, providing a finite-field analogue of (\ref{eq:Sha}):
\begin{equation}\label{eq:finiteSha}
\hat{I}_A(n)=\sum s_{\CO}(A)\hat{I}_{\CO}(n),\;\forall n\in\Lh^{nil}(k).
\end{equation}
Here the coefficients $s_{\CO}(A)$ are determined by the multiplicities and stalks of the aforementioned IC extensions. They are essentially what are called {\bf Green functions} \cite[\S5]{Spr76} and are known to be algorithmically computable in terms of characters of the Weyl group, at least when $\mathrm{char}(k)$ is good, see e.g. \cite[\S5]{Sho87}.

When transitioning from the finite-field realm to the $p$-adic realm, the Moy-Prasad filtration suggests that more than just the direct finite-field analogue is needed. It is shown in \cite[Thm. 4.1]{RY14} that a successive quotient of the Moy-Prasad filtration can be the algebraic representation $H^{\theta}\curvearrowright\Lh^{\theta=\lambda}$, where $\theta$ is a finite-order automorphism on a reductive group $H$, $\lambda$ is some eigenvalue of $\theta$ acting on $\Lh$ and the action is a restriction of adjoint representation. When $\theta$ is trivial and $\lambda=1$, we recover the adjoint representation $H\curvearrowright\Lh$ itself, which arises from integral-depth Moy-Prasad quotient, such as the depth-zero quotient. We refer to $H\curvearrowright\Lh$ as the ungraded case.

In general, we would like (\ref{eq:finiteSha}) in the {\bf graded case}, in which  $A\in\Lh^{\theta=\lambda}(k)$, $\CO\subset\Lh^{\theta=\lambda}(k)$ is an $H(k)$-orbit, $n$ is taken in $\Lh^{nil}(k)\cap\Lh^{\theta=\lambda^{-1}}$, and $I_A$ (resp. $I_{\CO}$) are sums over $H(k)$-conjugates of $A$ (resp. $H(k)$-orbit $\CO$). In the language of Springer theory, this means that we will get involved with $H^{\theta}$-equivariant perverse sheaves on $(\Lh^{\theta=\lambda})^{nil}$ as well as their Deligne-Fourier transforms, the latter known as {\it character sheaves}. The study of the ungraded case - named $H$-equivariant character sheaves on $\Lh$ - is called generalized Springer theory and is largely developed in \cite{Lus84}. We note that a majority of known wave-front set examples for depth-zero representations are based on Lusztig's work \cite{Lus92} which in turn relies on his renowned work on character sheaves \cite{Lus85} for the ``ungraded'' group $H\curvearrowright H$.

Character sheaves in the graded case can exhibit significantly different behaviors (see e.g. \cite{CVX18} and \cite{VX22}):
\begin{enumerate}[label=(\roman*)]
	\item The number of cuspidal objects can grow faster than any power of the rank, while in the ungraded case there is at most one per central character.
	\item Cuspidal character sheaves can be supported on the whole $\Lh^{\theta=\lambda}$, while in the ungraded case they are always supported on nilpotent mod center elements.
	\item These character sheaves are typically IC extensions of local systems of infinite order, whereas in the ungraded case, they always have finite order.
\end{enumerate}
A consequence of these new phenomena is that in the graded case, it is no longer possible to write the coefficients $s_{\CO}(A)$ in (\ref{eq:finiteSha}) in terms of discrete invariants such as characters of Weyl groups. Instead, the coefficients can depend on ``continuous invariants'' like point-counts for a family of varieties over $k$, such as the family $X_A^{(\ref{hard})}$ in (\ref{Sym2cover}) in the proof of Lemma \ref{hard}.

The graded case relevant to this paper is given by $H=GL_n$ and $\theta(g)=(g^t)^{-1}$, so that $H^{\theta}\curvearrowright\Lh^{\theta=\lambda}$ is precisely $O_n\curvearrowright\Sym^2(k^n)$, and this is the exact grading studied in \cite{CVX18}. In fact, it is possible to show that for $n=7$ and $A\in\Sym^2(k^7)$ in the introduction, the wave-front set for $(s_{\CO}(A))$ in the graded case of (\ref{eq:finiteSha}) agrees with Theorem \ref{main}. This is at least heuristically a necessary condition for Theorem \ref{main}. 

In general, the theory of character sheaves in the graded case is still being developed and is far less complete than the ungraded case. In particular, we don't fully understand the aforementioned ``continuous invariants'' yet. However, knowing that $|X_A^{(\ref{hard})}(k)|$ will appear and having examples about such appearance allow us to study something. Since character sheaves in the graded case can look quite different from the ungraded case, we would say whatever is possibly computable in terms of them (such as wave-front set for rational-depth representations) can, in principle, behave very differently from ungraded case (or the wave-front set for depth-zero representations). For example, we are inclined to expect the wave-front set of any depth-zero representation to be contained in a single geometric orbit.

%We believe that the different behavior in graded Springer theory implies that harmonic analysis for rational depths can behave differently from the situation for integral depths, in a way similar to the phenomenon in \cite{VX18a} that awaits to be made precise. In particular the local character expansion of a rational-depth representation can behave differently from those of depth-zero representations. Meanwhile, it seems still possible that the geometric wave-front set of any depth-zero representation is a singleton, as shown in many cases in \cite{Wal18}, \cite{Wal20}, \cite{CMBO21}, \cite{AGS22a} and \cite{CMBO22a}.

We remark that stalks of character sheaves on our grading $O_n\curvearrowright\Sym^2(k^n)$ probably first implicitly appear in Hales' work \cite{Ha94} where he discovered that point-counts on hyperelliptic curves appear in the stable subregular Shalika germs. Later my advisor B. Gross noticed explicitly that this grading should lead to interesting results and told me about it in 2012. Our work is intellectually in debt to them.

\subsection{Arithmetic geometry}\label{subsec:AG} Point-counts on varieties over the residue field have long been observed to appear in $p$-adic integrals. For the coefficients in (\ref{eq:Sha}) as well as the graded case of (\ref{eq:finiteSha}), when $p \gg \rank H$ these varieties can be chosen from a family of Hessenberg varieties \cite[\S2.5]{GKM06} or their \'{e}tale neighborhoods in \cite[\S4]{Ts15d}. Hessenberg varieties have proven to be useful in graded and affine Springer theory, as demonstrated in \cite{GKM06}, \cite{OY16}, \cite{LY17}, and \cite{CVX18}.

The variety $X^{(\ref{hard})}_A$ in (\ref{Sym2cover}) is an example of a Hessenberg variety, which is a projective smooth surface \cite[\S2.5]{GKM06} and in fact also geometrically connected. More importantly, as a Hessenberg variety it has the form:
\begin{equation}\label{Hessenbergform}
X_A(k):=\{g\in SO_7(k)/B(k)\;|\;g^{-1}Ag\text{ lives in some subspace of }\Sym^2(k^7)\}.
\end{equation}
We see that $X_A$ has a natural left $\operatorname{Stab}_{SO_7}(A)$-action.
%Meanwhile, recall that we took $\til{A}$ any lift in $\Lg(F)_{-\frac{1}{2}}$ of $A$ under (\ref{map}). It happens in our case that $T:=Z_G(\til{A})$ is anisotropic and the group scheme $X_*(T)_{I_F}$ is canonically isomorphic \cite[proof of Lemma 4.5]{Ts17} to $\Stab_{O_7}(A)$ in Lemma \ref{stabilizer}. We have $\ker(H^1(F,T)\ra H^1(F,G))=\ker(H^1(k,\Stab_{O_7}(A))\ra H^1(k,O_7))=H^1(k,\Stab_{SO_7}(A))$.
Let $C:=\Stab_{SO_7}(A)$. Each class $\alpha\in H^1(k,C)$ parameterizes an $O_7(k)$-conjugacy class of $A_{\alpha}\in\Sym^2(k^7)$ that is conjugate to $A$ under $O_7(\bar{k})$.
%, and a $G(F)$-conjugacy class in $\til{A}_{\alpha}\in\Lg(F)$ that is conjugate to $A'$ under $G(F^{sep})$. It is in fact possible to arrange that $\til{A}_{\alpha}$ is again a lift of $A_{\alpha}$ in $\Lg(F)_{-\frac{1}{2}}$.
From (\ref{Hessenbergform}) one can deduce that
\[
X_{A_{\alpha}}\cong X_A\times^C\alpha:=C\bsl (X_A\times\alpha),
\]
where $\alpha\in H^1(k,C)$ is realized as a $C$-torsor over $\Spec k$ and $C$ acts on $X_A\times\alpha$ via the diagonal action. Recall that in Lemma \ref{stabilizer} we saw $C$ is an abelian $2$-torsion group scheme of order $2^{7-1}=64$. Consider the hypotheses:
\begin{enumerate}[label=(\roman*)]
	\item\label{assump:split} $C$ is a constant group scheme, so that $|C(k)|=|H^1(k,C)|=2^{7-1}$.
	\item\label{assump:free} $C$ acts freely on $X_A$.
	\item\label{assump:less} The quotient $k$-variety $C\bsl X_A$ has fewer than $2^{7-1}$ rational points.
\end{enumerate}
Among these, only hypothesis \ref{assump:split} is straightforward to confirm: by Lemma \ref{stabilizer}, it is equivalent to the characteristic polynomial of $A$ splitting completely in $k$. Assuming hypothesis \ref{assump:free}, each fiber of $X_A\ra C\bsl X_A$ above a $k$-point is a $C$-torsor over $\Spec k$. Replacing $X_A$ by $X_A\times^C\alpha$ corresponds to twisting every such fiber by $\alpha$. Given hypotheses \ref{assump:split} and \ref{assump:less}, there must be some choice of $\alpha$ for which all rational fibers of $X_A\ra C\bsl X_A$ are non-trivial torsors, i.e. none of them has a rational point. This implies that $X_{A_{\alpha}}(k)=\emptyset$, which is what we need in Lemma \ref{hard} and serves as the most important input for our construction.

The hypotheses themselves can't always be made true, but in \cite{Ts17} these hypotheses are verified for some double cover of $X_A^{(\ref{hard})}$ called\footnote{It's defined before \cite[Thm. 3.4]{Ts17} and shown to be highly related to Hales' hyperelliptic curves. The argument after it is pretty long, though.} $F_{T,2}$, in fact for $U_7$ replaced by any ramified $U_{2n+1}$, in which case $X_A^{(\ref{hard})}$ is replaced by the closed subvariety of $SO_{2n+1}/B$ parametrizing $gB$ for which $g^{-1}Ag$ has zero entries below the sub-diagonal except for the middle $3$ entries right below the sub-diagonal, similar to Lemma \ref{hard}. Based on the arguments above, this gives the following reinterpretation of \cite[Cor. 6.3]{Ts17}:\hp

\stepcounter{equation}
\newcounter{EquationCounter}
\setcounter{EquationCounter}{\value{equation}} % Store the current equation counter value

\begin{enumerate}[leftmargin=.85cm]
	\item[(\arabic{EquationCounter})]\label{Ts17} 
	 {\it For any integer $d>0$, there exists an integer $N(d)$ such that for any $n\ge N(d)$ and any odd prime number $q\le n^d$, there exists a regular semisimple $A'\in \Sym^2(\mathbb{F}_q^{2n+1})$ such that the wave-front sets of the supercuspidal representations of $U_{2n+1}(E/F)$ produced from $A'$ excludes nilpotent orbits of Jordan type $(n)$, $(n-1,1)$ and $(n-2,2)$.}
\end{enumerate}\qp

Here, $\Sym^2(\mathbb{F}_q^{2n+1})$ is the same as $\Sym^2(k^7)$ in the introduction, except that $7$ is replaced by $2n+1$, $k=\mathbb{F}_3$ is replaced by $k=\mathbb{F}_q$, and $F=\Q_3$ is replaced by any $p$-adic field with residue field $k$. The construction of the supercuspidal representations is the same as in the introduction. We emphasize that the ingredient hypothesis \ref{assump:less} is valid for large enough $n$ and $q\le n^d$ because the order of $C$ grows exponentially in $n$ while $|(C\bsl X_{A'})(k)|=O(q^2n^2)$.

Given the existence of such supercuspidals, our heuristic is then that ``almost always'' their wave-front sets will include some nilpotent orbits of Jordan type $(n-2,1,1)$ and $(n-3,3)$, and in particular not all in a single geometric orbit. In light of Lemma \ref{triv} and \ref{Jac} (and the way they are used in Proposition \ref{ss}) this is asking for some other varieties $X_{A'}^{(\ref{triv})}$ and $X_{A'}^{(\ref{Jac})}$ over $k$ to have a $k$-point. Unfortunately, while we feel like this holds for almost all random samples, we lack a conceptual proof confirming its validity for at least one sample. As a result, to generate a (counter-)example, we rely on computer programs, which successfully produce one already for the minimal case $n=3$ and $q=3$.

Because of this heuristic, we expect similar counterexamples also in split type {\bf B}, {\bf C}, {\bf D} - with arbitrarily large ranks and suitable values of $q$, similar to those in (\arabic{EquationCounter}).

\p
\bibliographystyle{amsalpha}
\bibliography{biblio.bib}
\p

\appendix
\section{Computer programs}\label{AppReadme}
Two programs are included in this appendix. Both programs are in {\bf Magma} and can be run by copying the code onto \url{http://magma.maths.usyd.edu.au/calc/}. The second program is almost self-explained, while the first program needs a bit of explanation as to why it proves Lemma \ref{hard}. In the proof of Lemma \ref{hard} we had to prove that the variety
\begin{equation}\tag{\ref{Sym2cover}}\label{Sym2covercopy}
	X^{(\ref{hard})}_A(k):=\{g\in SO_7(k)/B(k)\;|\;g^{-1}Ag\in
	\matr{
		*&*&*&*&*&*&*\\
		*&*&*&*&*&*&*\\
		0&*&*&*&*&*&*\\
		0&*&*&*&*&*&*\\
		0&0&*&*&*&*&*\\
		0&0&0&*&*&*&*\\
		0&0&0&0&0&*&*
	}
	\}.
\end{equation}
has no rational point. Here $SO_7(k)$ is defined with respect to the form
\begin{equation}\tag{\ref{form2}}\label{form2copy}
	\langle x,y\rangle=x_1y_7+...+x_7y_1,
\end{equation}
and $B\subset SO_7$ is the Borel subgroup of upper triangular matrices. The set $SO_7(k)/B(k)=(SO_7/B)(k)$ parameterizes flags $0\subset V_1\subset V_2\subset V_3\subset V_4\subset V_5\subset V_6\subset k^7$ such that $\dim_kV_i=i$ and that
\begin{equation}\tag{i}\label{eq:i}
	\langle V_i,V_j\rangle =0,\;\forall i,j\text{ such that }i+j\le 7.
\end{equation}
More precisely, a coset $g\in SO_7(k)/B(k)$ parametrizes the flag given by $(V_i=\mathrm{span}(ge_1,...,ge_i))$, where $e_i$ is the column vector with $1$ on the $i$-th entry and $0$ elsewhere.

The condition that $g^{-1}Ag$ is of the form indicated on the right of (\ref{Sym2covercopy}) is, by unfolding the definition, equivalent to that
\[
\begin{array}{c}
	A(V_1)\subset V_2,\\
	A(V_2)\subset V_4,\\
	A(V_3)\subset V_5,\\
	A(V_5)\subset V_6.
\end{array}
\]
Here $V_i$ are as above, e.g. $V_1=\mathrm{span}(ge_1)$. These conditions imply $A(ge_1)\in V_2$ and $A^2(ge_1)\in V_4$. Thanks to (\ref{eq:i}) we have $\langle ge_1,ge_1\rangle=\langle ge_1,Age_1\rangle=\langle Age_1,Age_1\rangle=\langle Age_1,A^2ge_1\rangle=0$. Moreover, $A^2(ge_1)\in V_4=\mathrm{span}(ge_3)+V_3\implies A^2ge_1\in cge_4+V_3$ for some $c\in k$. Since $\langle ge_4,ge_4\rangle=\langle e_4,e_4\rangle=1$, this implies $\langle A^2ge_1,A^2ge_1\rangle=c^2$ is a square in $k$. This gives the following:\vskip.2cm

\noindent{\bf Lemma.} A sufficient condition for $X^{(\ref{hard})}_A(k)=\emptyset$ is that there does not exist non-zero $v\in k^7$ such that $\langle v,v\rangle=\langle v,Av\rangle=\langle Av,Av\rangle=\langle Av,A^2v\rangle=0$ and $\langle A^2v,A^2v\rangle\in k$ is a square.\vskip.2cm

It is possible to prove that this is also a necessary condition. In any case, this is an algorithm that's easy to program, which we do in \ref{App2.2} and briefly explain below:

\includegraphics{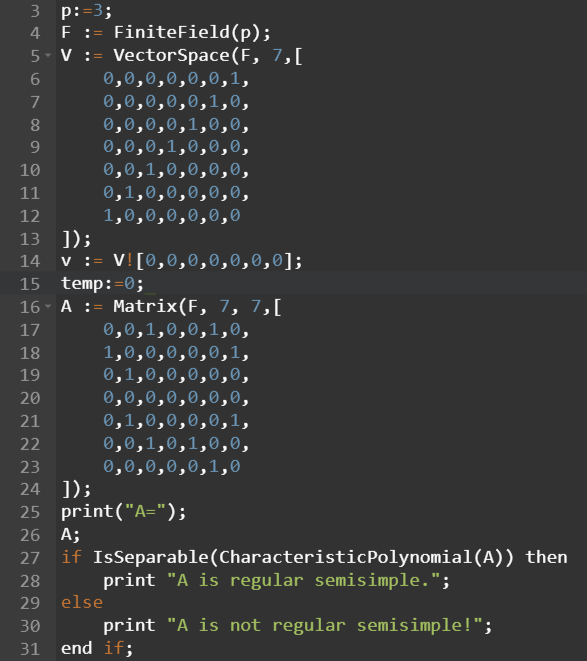}
In this part we set up V to be the our inner product space $k^7$ with the inner product given by form (\ref{form2copy}) and set up $v$ to be a vector in $V$. Next we set up matrix $A$ and verify that $A$ has separable characteristic polynomial. Lastly we print all these results for our convenience.

\includegraphics{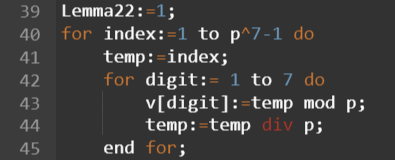}We begin by the auxiliary variable Lemma22 (``Lemma 2.2 is assumed true unless we find such a $v$.'') Then, we enumerate all $v\in k^7\backslash\{0\}$ by running over all $3^7-1$ non-zero vectors. The variable {\bf index} is made to run from $1$ to $3^7-1$, and the vector-variable {\bf v[-]}, namely $v$, is the vector given by the $3$-adic expression of {\bf index}.

\includegraphics{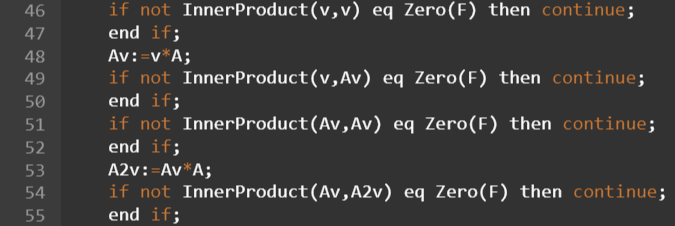}This part says ``skip to the next {\bf index} (i.e. next possible $v$) whenever any of $\langle v,v\rangle$, $\langle v,Av\rangle$, $\langle Av,Av\rangle$ and $\langle Av,A^2v\rangle$ is non-zero. The program uses row vectors instead of column vectors, so $A$ is multiplied from the right. Note that $A$ is self-adjoint so row vectors and column vectors can be identified, under the transpose map which sends $\matr{1\\2\\3}$ to $\matr{3&2&1}$.

\includegraphics{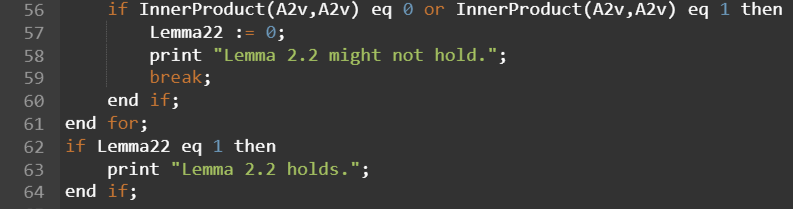}Lastly, it checks if $\langle A^2v,A^2v\rangle$ is a square (in $\mathbb{F}_3$ the squares are $0$ and $1$). If it is, then it says that the sufficient condition is not met and halts the program (``break'' the loop). If this never happens and the loop for {\bf index} is run til the end (namely, all $v\in k^7\backslash\{0\}$ are verified), it outputs that Lemma \ref{hard} holds, because our sufficient condition has been verified.

\includegraphics{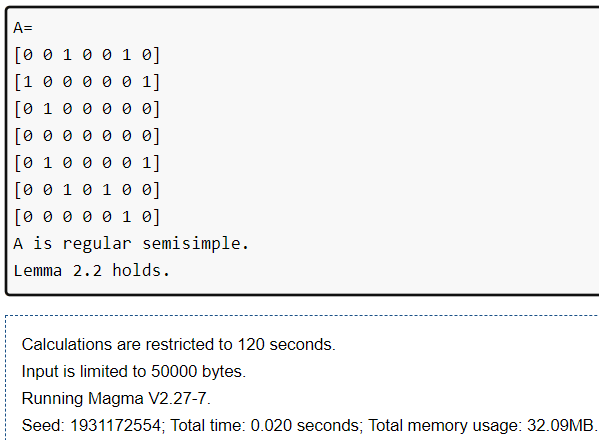} The output on the Magma website is, as asserted, affirmative.\vskip.2cm

For $p=5$, we also have the same wave-front set for a supercuspidal representation of $U_7(E/F)$ with $F=\mathbb{Q}_5$ given by
\[
A=\matr{
	0&0&1&1&1&0&0\\
	1&0&0&0&0&0&0\\
	0&1&0&0&0&0&1\\
	0&0&0&0&0&0&1\\
	0&1&0&0&0&0&1\\
	0&0&1&0&1&0&0\\
	0&0&0&0&0&1&0
}\in\Sym^2(\mathbb{F}_5^7)
\]
which can be verified using the same program (changing $p:=5$ and changing the square verification). Similar for $p:=7$ and
\[
A=\matr{
	0&0&4&2&2&0&0\\
	1&0&0&0&0&0&0\\
	0&3&0&0&0&0&2\\
	0&0&0&0&0&0&2\\
	0&1&0&0&0&0&4\\
	0&0&1&0&3&0&0\\
	0&0&0&0&0&1&0
}\in\Sym^2(\mathbb{F}_7^7).
\]
We can't find any example for $p=11$. We'd guess there is none.

\subsection{Magma code for Lemma \ref{hard}}\label{App2.2}
\small
\begin{Verbatim}[tabsize=4, commandchars=\\\{\}]
p:=\textcolor{purple}{3};
F := FiniteField(p);
V := VectorSpace(F, \textcolor{purple}{7},[
	\textcolor{purple}{0,0,0,0,0,0,1,}
	\textcolor{purple}{0,0,0,0,0,1,0,}
	\textcolor{purple}{0,0,0,0,1,0,0,}
	\textcolor{purple}{0,0,0,1,0,0,0,}
	\textcolor{purple}{0,0,1,0,0,0,0,}
	\textcolor{purple}{0,1,0,0,0,0,0,}
	\textcolor{purple}{1,0,0,0,0,0,0}
]);
v := V![\textcolor{purple}{0,0,0,0,0,0,0}];
temp:=0;
A := Matrix(F, \textcolor{purple}{7}, \textcolor{purple}{7},[
	\textcolor{purple}{0,0,1,0,0,1,0,}
	\textcolor{purple}{1,0,0,0,0,0,1,}
	\textcolor{purple}{0,1,0,0,0,0,0,}
	\textcolor{purple}{0,0,0,0,0,0,0,}
	\textcolor{purple}{0,1,0,0,0,0,1,}
	\textcolor{purple}{0,0,1,0,1,0,0,}
	\textcolor{purple}{0,0,0,0,0,1,0}
]);
print(\textcolor{teal}{"A="});
A;
\textcolor{brown}{if} IsSeparable(CharacteristicPolynomial(A)) \textcolor{brown}{then}
	print \textcolor{teal}{"A is regular semisimple."};
\textcolor{brown}{else}
	print \textcolor{teal}{"A is not regular semisimple!"};
\textcolor{brown}{end if};

Lemma22:=\textcolor{purple}{1};
\textcolor{brown}{for} index:=\textcolor{purple}{1} to p^\textcolor{purple}{7}-\textcolor{purple}{1} \textcolor{brown}{do}
	temp:=index;
	\textcolor{brown}{for} digit:= \textcolor{purple}{1} to \textcolor{purple}{7} \textcolor{brown}{do}
		v[digit]:=temp mod p;
		temp:=temp div p;
	\textcolor{brown}{end for};
	\textcolor{brown}{if} not InnerProduct(v,v) eq \textcolor{purple}{0} \textcolor{brown}{then continue};
	\textcolor{brown}{end if};
	Av:=v*A;
	\textcolor{brown}{if} not InnerProduct(v,Av) eq \textcolor{purple}{0} \textcolor{brown}{then continue};
	\textcolor{brown}{end if};
	\textcolor{brown}{if} not InnerProduct(Av,Av) eq \textcolor{purple}{0} \textcolor{brown}{then continue};
	\textcolor{brown}{end if};
	A2v:=Av*A;
	\textcolor{brown}{if} not InnerProduct(Av,A2v) eq \textcolor{purple}{0} \textcolor{brown}{then continue};
	\textcolor{brown}{end if};
	\textcolor{brown}{if} InnerProduct(A2v,A2v) eq \textcolor{purple}{0} \textcolor{brown}{or} InnerProduct(A2v,A2v) eq \textcolor{purple}{1} \textcolor{brown}{then}
		Lemma22 := \textcolor{purple}{0};
		print \textcolor{teal}{"Lemma 2.2 might not hold."};
		break;
	\textcolor{brown}{end if};
\textcolor{brown}{end for};
\textcolor{brown}{if} Lemma22 eq \textcolor{purple}{1} \textcolor{brown}{then}
	print \textcolor{teal}{"Lemma 2.2 holds."};
\textcolor{brown}{end if};
\end{Verbatim}

\subsection{Magma code for Lemma \ref{Jac}}\label{App2.4}
\begin{Verbatim}[tabsize=4, commandchars=\\\{\}]
p:=3;
F := FiniteField(p);
A := Matrix(F, \textcolor{purple}{7}, \textcolor{purple}{7},[
	\textcolor{purple}{0,0,1,0,0,1,0,}
	\textcolor{purple}{1,0,0,0,0,0,1,}
	\textcolor{purple}{0,1,0,0,0,0,0,}
	\textcolor{purple}{0,0,0,0,0,0,0,}
	\textcolor{purple}{0,1,0,0,0,0,1,}
	\textcolor{purple}{0,0,1,0,1,0,0,}
	\textcolor{purple}{0,0,0,0,0,1,0}
]);
h := Matrix(F, \textcolor{purple}{7}, \textcolor{purple}{7},[
	\textcolor{purple}{2,1,1,2,1,2,2,}
	\textcolor{purple}{1,1,0,0,0,1,2,}
	\textcolor{purple}{1,0,1,2,0,0,1,}
	\textcolor{purple}{1,2,2,0,1,1,1,}
	\textcolor{purple}{0,2,0,0,0,0,1,}
	\textcolor{purple}{1,1,2,1,0,1,0,}
	\textcolor{purple}{0,0,1,1,1,0,0}
]);
g:=ZeroMatrix(F,\textcolor{purple}{7},\textcolor{purple}{7});
\textcolor{brown}{for} i:=\textcolor{purple}{1} to \textcolor{purple}{7} \textcolor{brown}{do}
	\textcolor{brown}{for} j:=\textcolor{purple}{1} to \textcolor{purple}{7} \textcolor{brown}{do}
		g[i][j]:=h[\textcolor{purple}{8}-j][\textcolor{purple}{8}-i];
	\textcolor{brown}{end for};
\textcolor{brown}{end for};
print(\textcolor{teal}{"g="});
g;
print(\textcolor{teal}{"h*g="});
h*g;
print(\textcolor{teal}{"h*A*g="});
h*A*g;
\end{Verbatim}
\p

\end{document}